\documentclass[letterpaper, 10 pt, conference]{ieeeconf}  

\IEEEoverridecommandlockouts                              
\usepackage{amsmath,times,amssymb}
\usepackage{cite}
\usepackage{psfrag,graphicx}
\usepackage{changebar}

\allowdisplaybreaks

\newtheorem{theorem}{Theorem}[section]
\newtheorem{lemma}[theorem]{Lemma}
\newtheorem{algorithm}[theorem]{Algorithm}
\newtheorem{assumption}[theorem]{Assumption}
\newtheorem{proposition}[theorem]{Proposition}
\newtheorem{definition}[theorem]{Definition}

\newtheorem{remark}[theorem]{Remark}

\label{newcom}
\newcommand{\xb}{\textbf{x}}
\newcommand{\ub}{\textbf{u}}

\newcommand{\Xc}{\mathcal{X}}
\newcommand{\Uc}{\mathcal{U}}
\newcommand{\Xf}{\mathcal{X}_\mathrm{f}}

\newcommand{\rr}{\mathbb{R}}
\newcommand{\zz}{\mathbb{Z}}

\newcommand{\Nr}{\mathcal{N}}

\newcommand{\grad}{\nabla}

\newcommand{\Lc}{\mathcal{L}}

\newcommand{\cvd}{\hfill $\Box$}

\newcommand{\Omegab}{\mathbf{\Omega}}

\title{{\LARGE \bf
A distributed optimization-based approach for \\ hierarchical model predictive control of large-scale systems with \\ coupled dynamics and constraints}\\
{\small (extended version of the CDC-ECC'11 paper, with proofs)}
}

\author{Minh Dang Doan, Tam\'{a}s Keviczky, and Bart De Schutter
\thanks{The authors are with Delft Center for Systems and Control, Delft University of Technology,
   Delft, The Netherlands
        {\tt\small \{m.d.doan, t.keviczky, b.deschutter\}@tudelft.nl}}%
}

\begin{document}

\maketitle

\thispagestyle{empty}
\pagestyle{empty}

\begin{abstract}

We present a hierarchical model predictive control approach for large-scale systems based on dual decomposition. The proposed scheme allows coupling in both dynamics and constraints between the subsystems and generates a primal feasible solution within a finite number of iterations, using primal averaging and a constraint tightening approach. The primal update is performed in a distributed way and does not require exact solutions, while the dual problem uses an approximate subgradient method. Stability of the scheme is established using bounded suboptimality.

\end{abstract}

\section{Introduction}

Coordination and control of interacting subsystems is an essential requirement for optimal operation and enforcement of critical operational constraints in large-scale industrial processes and infrastructure systems \cite{rawlings_coordinating_2008}. Model Predictive Control (MPC) has become the method of choice when designing control systems for such applications \cite{Mac:02, CamBor:99, RawMay:09}, due to its ability to handle important process constraints explicitly. MPC relies on solving finite-time optimal control problems repeatedly online, which may become prohibitive for large-scale systems due to the problem size or communication constraints. Recent efforts have been focusing on how to decompose the underlying optimization problem in order to arrive at a distributed or hierarchical control system that can be implemented under the prescribed computational and communication limitations \cite{scattolini_architecture_2009, Venkat_tcst:2008}. One common way to decompose an MPC problem with coupled dynamics or constraints is to use dual decomposition methods \cite{yuji_wakasa_decentralized_2008, NecSuy_prox:2008, DoaKev:09}, which typically lead to iterative algorithms (in either a distributed or a hierarchical framework) that converge to feasible solutions only asymptotically. Implementing such approaches within each MPC update period can be problematic for some applications.

Recently, we have presented a dual decomposition scheme for solving large-scale MPC problems with coupling in both dynamics and constraints, where primal feasible solutions can be obtained even after a finite number of iterations~\cite{DoaKev:11IFACWC}. In the current paper we present a novel method that is motivated by the use of constraint tightening in robust MPC~\cite{Kuwata:07}, along with a primal averaging scheme and distributed Jacobi optimization. Since an exact optimum of the Lagrangian is not assumed to be computable in finitely many iterations, an approximate scheme is needed for solving the MPC optimization problem at each time step.
We present a solution approach that requires a nested two-layer iteration structure and the sharing of a few crucial parameters in a hierarchical fashion.
The proposed framework guarantees primal feasible solutions and MPC stability using a finite number of iterations with bounded suboptimality.

The paper is organized as follows. In Section \ref{sec_problem}, we describe the MPC optimization problem and its tightened version, which will be used to guarantee feasibility of the original problem even with a suboptimal primal solution. Section~\ref{sec_algorithm} describes the main elements of the algorithm used to solve the dual version of the tightened optimization problem: the approximate subgradient method and the distributed Jacobi updates. In Section~\ref{sec_properties}, we show that the primal average solution generated by the approximate subgradient algorithm is a feasible solution of the original optimization problem, and that the cost function decreases through the MPC updates. This allows it to be used as a Lyapunov function for showing closed-loop MPC stability. Section~\ref{sec_conclusion} concludes the paper and outlines future research.

\section{Problem description}\label{sec_problem}
\subsection{MPC problem}
We consider $M$ interconnected subsystems with coupled discrete-time linear time-invariant dynamics:
\begin{align} \label{eq_dmodel}
x^i_{k+1} &= \sum_{j=1}^M A^{ij} x^j_k + B^{ij} u^j_k,  \quad i=1,\dots,M 
\end{align}
and the corresponding centralized state-space model:
\begin{align} \label{eq_cmodel}
x_{k+1} &= A x_k + B u_k
\end{align}
with $x_k=[(x_k^1)^T (x_k^2)^T \dots (x_k^M)^T]^T, u_k=[(u_k^1)^T (u_k^2)^T \dots (u_k^M)^T]^T$, $A=[A_{ij}]_{i,j \in \{1,\dots,M\}}$ and $B=[B_{ij}]_{i,j \in \{1,\dots,M\}}$.

The MPC problem at time step $t$ is formed using a convex cost function and convex constraints:

\begin{align}
\min_{\ub, \xb} \quad & \sum_{k=t}^{t+N-1} \bigg( x_k^T Q x_k + u_k^T R u_k \bigg) + x_{t+N}^T P x_{t+N} \label{eq_cmpc}\\
 \textrm{s.t.} \quad & x^i_{k+1} = \sum_{j \in \Nr^i} A^{ij} x^j_k + B^{ij} u^j_k, \nonumber \\
&\quad\quad\quad\quad i=1,\dots,M, \quad k=t,\dots,t+N-1 \label{eq_cmpc_model_const} \\
& x_k \in \Xc, k=t+1,\dots,t+N-1 \label{eq_cmpc_x_const} \\
& x_{t+N} \in \Xf \subset \Xc \label{eq_cmpc_xf_const} \\
& u_k \in \Uc, k=t,\dots,t+N-1 \label{eq_cmpc_u_const} \\
& u^i_k \in \Omega_i, i=1,\dots,M, \quad k=t,\dots,t+N-1 \label{eq_cmpc_u_local_const} \\
& x_t = x(t) \in \Xc \label{eq_cmpc_xt}
\end{align}
where $\ub=[u_t^T,\dots,u_{t+N-1}^T]^T$, $\xb=[x_{t+1}^T,\dots,x_{t+N}^T]^T$, the matrices $Q$, $P$, and $R$ are block-diagonal and positive definite, the constraint sets $\Uc$, $\Xc$ and $\Xf$ are polytopes and have nonempty interiors, and each local constraint set $\Omega_i$ is a hyperbox. Each subsystem $i$ is assigned a neighborhood, denoted $\Nr^i$, containing subsystems that have direct dynamical interactions with subsystem $i$, including itself. The initial state $x_t$ is the current state at time step $t$.

As $\Uc$, $\Xc$ and $\Xf$ are polytopes, the constraints \eqref{eq_cmpc_x_const} and \eqref{eq_cmpc_xf_const} are represented by linear inequalities. Moreover, the state vector $\xb$ is affinely dependent on $\ub$. Hence, we can eliminate state variables $x_{t+1},\dots,x_{t+N}$ and transform the constraints \eqref{eq_cmpc_model_const}, \eqref{eq_cmpc_x_const}, and \eqref{eq_cmpc_xf_const} into linear inequalities of the input variable $\ub$. Eliminating the state variables in \eqref{eq_cmpc}--\eqref{eq_cmpc_xt} leads to an optimization problem in the following form:
\begin{align} \label{eq_opt0}
 f_t^* = \min_{\ub} \quad & f(\ub, x_t) \\
 \textrm{s.t.} \quad & g(\ub, x_t) \leq 0 \label{eq_opt0_const} \\
& \ub \in \Omegab \label{eq_opt0_domain}
\end{align}
where $f$ and $g=[g_{1}, \dots, g_{m}]^T$ are convex functions, and $\Omegab=\prod_{i=1}^M \Omegab_i$ with each $\Omegab_i = \prod_{k=0}^{N-1} \Omega_i$ is a hyperbox. Note that $f(\ub, x_t) > 0, \forall \ub \neq 0, x_t \neq 0$, due to the positive definiteness of $Q$, $P$, and $R$.

We will use $(\ub_{t}, \xb_{t})$ to denote a feasible solution generated by the controller for problem \eqref{eq_cmpc}--\eqref{eq_cmpc_xt} at time step $t$. This solution is required to be feasible but not necessarily optimal.
We will make use of the following assumptions:

\begin{assumption}\label{asm_decentralized_controller}
There exists a block-diagonal feedback gain $K$ such that the matrix $A+BK$ is Schur (i.e., a decentralized stabilizing control law for the unconstrained aggregate system).
\end{assumption}

\begin{assumption}\label{asm_pos_inv_terminal_set}
The terminal constraint set $\Xf$ is positively invariant for the closed-loop $x_{k+1}=(A+BK)x_k$ ($x \in \mathrm{int}(\Xf) \Rightarrow (A+BK)x \in \mathrm{int}(\Xf)$).
\end{assumption}

\begin{assumption}\label{asm_cond_slater}
The Slater condition holds for problem \eqref{eq_opt0}--\eqref{eq_opt0_domain}, i.e., there exists a vector that satisfies strict inequality constraints \cite{Ber_nlp_book:1999}. It is also assumed that prior to each time step $t$, a Slater vector $\bar{\ub}_t$ is available, such that
\begin{align} \label{eq_cond_slater}
 g_{j}(\bar{\ub}_t,x_t) < 0, j=1,\dots,m
\end{align}
\end{assumption}

\begin{remark}
Since $g(\ub, x_t) \leq 0$ has a nonempty interior, so do its components $g_{j}(\ub, x_t) \leq 0, j=1,\dots,m$. Hence, there will always be a vector that satisfies the Slater condition \eqref{eq_cond_slater}. In fact, we will only need to find the Slater vector $\bar{\ub}_0$ for the first time step, which can be computed off-line. In Section~\ref{sec_slate_update} we will show that a new Slater vector can then be obtained for each $t \geq 1$, using Assumption~\ref{asm_pos_inv_terminal_set}.
\end{remark}

\begin{assumption}\label{asm_cost_decrease}
At each time step $t$, the following holds
  \begin{align}\label{eq_cost_decrease}
   f(\ub_{t-1},x_{t-1}) - f(\bar{\ub}_{t},x_{t}) > x_{t-1}^T Q x_{t-1} + u_{t-1}^T R u_{t-1}
  \end{align}

For later reference, we define $\Delta_t > 0$ which can be computed before time step $t$ as follows:
\begin{align}\label{eq_def_Delta}
\Delta_t = x_{t-1}^T Q x_{t-1} + u_{t-1}^T R u_{t-1}
\end{align}
\end{assumption}

\begin{remark}
Assumption~\ref{asm_cost_decrease} is often satisfied with an appropriate terminal penalty matrix $P$. A method to construct a block-diagonal $P$ with a given decentralized stabilizing control law is provided in \cite{Siljak:78_lss}.
\end{remark}

\begin{assumption}\label{asm_g_norm_bound}
For each $x_t \in \Xc$, the Euclidean norm of $g(\ub,x_t)$ is bounded:
\begin{align}\label{eq_cond_bound}
 L_t \geq \|g(\ub,x_t)\|_2, \forall \ub \in \Omegab
\end{align}
\end{assumption}

\begin{remark}
In the first time step, with given $x_0$, we can find $L_0$ by evaluating $\|g(\ub,x_0)\|_2$ at the vertices of $\Omegab$, the maximum will then satisfy \eqref{eq_cond_bound} for $t=0$, due to the convexity of $g$ and $\Omegab$. For the subsequent time steps, we will present a simple method to update $L_t$ in Section~\ref{sec_norm_bound_update}.
\end{remark}

\subsection{The tightened problem}\label{sec_tightening}

We will not solve problem \eqref{eq_opt0}--\eqref{eq_opt0_domain} directly. Instead, we will make use of an iterative algorithm based on a tightened version of \eqref{eq_opt0}--\eqref{eq_opt0_domain}. Consider the tightened constraint:
\begin{align}\label{eq_tight_const}
 g'(\ub,x_t) \triangleq g(\ub,x_t) + \mathbf{1}_m c_t \leq 0
\end{align}
with $g'(\ub,x_t)=[g'_{1}, \dots, g'_{m}]^T$, $0 < c_t < \min_{j=1,\dots,m}\{-g_j(\bar{\ub}_t,x_t)\}$, and $\mathbf{1}_m$ the column vector with every entry equal to $1$. Due to \eqref{eq_cond_slater}, we have
 \begin{align} 
  \max_{j=1,\dots,m}\{g'_j(\bar{\ub}_t,x_t)\} = \max_{j=1,\dots,m}\{g_j(\bar{\ub}_t,x_t)\}+c_t < 0
 \end{align}

Hence $g'_{j}(\bar{\ub}_t,x_t) < 0, j=1,\dots,m$. Moreover, using \eqref{eq_cond_bound} and the triangle inequality of the 2-norm, we will get $L'_t = L_t + c_t$ as the norm bound for $g'$, i.e. $ L'_t \geq \|g'(\ub,x_t)\|_2, \forall \ub \in \Omegab$. Note that $L'_t$ implicitly depends on $x_t$, as $\bar{\ub}_t$ and $c_t$ are updated based on the current state $x_t$.

Using the tightened constraint \eqref{eq_tight_const}, we formulate the tightened problem:
\begin{align} \label{eq_opt1}
 {f'_t}^* = \min_{\ub} \quad & f(\ub, x_t) \\
 \textrm{s.t.} \quad & g'(\ub, x_t) \leq 0 \label{eq_opt1_const} \\
& \ub \in \Omegab \label{eq_opt1_domain}
\end{align}

\begin{remark}
Only the coupled constraints \eqref{eq_opt0_const} are tightened, while the local input constraints \eqref{eq_opt0_domain} are unchanged. The Slater condition also holds for the tightened problem \eqref{eq_opt1}--\eqref{eq_opt1_domain}, with $\bar{\ub}_t$ being the Slater vector. 
\end{remark}

\section{The proposed optimization algorithm}\label{sec_algorithm}

Our objective is to calculate a feasible solution for problem \eqref{eq_cmpc}--\eqref{eq_cmpc_xt} using a method that is favorable for distributed computation. The main idea is to use dual decomposition for the tightened problem \eqref{eq_opt1}--\eqref{eq_opt1_domain} instead of the original one, such that after a finite number of iterations the constraint violations in the tightened problem will be less than the difference between the tightened and the
original constraints. Thus, even after a finite number of iterations, we will obtain a primal feasible solution for the original MPC optimization problem.

\subsection{The dual problem}

We will tackle the dual problem of \eqref{eq_opt1}--\eqref{eq_opt1_domain}, in order to deal with coupled constraint $g'(\ub, x_t) \leq 0$ in a distributed way. In this section, we define the dual problem and its subgradient. For simplicity, in this section the dependence of functions on the initial condition $x_t$ is not indicated explicitly.

The Lagrangian of problem \eqref{eq_opt1}--\eqref{eq_opt1_domain} is defined as:
\begin{align}
 \mathcal{L}'(\ub, \mu) = f(\ub) + \mu^Tg'(\ub)
\end{align}
in which $\ub \in \Omegab, \mu \in \rr^m_+$.

The dual function for \eqref{eq_opt1}--\eqref{eq_opt1_domain}:
\begin{align}
 q'(\mu) = \min_{\ub \in \Omegab} \mathcal{L}'(\ub, \mu)
\end{align}
is a concave function on $\rr^m_+$, and it is non-smooth when $f$ and $g'$ are not strictly convex functions \cite{Ber_nlp_book:1999}.

Given the assumption that Slater condition holds for \eqref{eq_opt1}--\eqref{eq_opt1_domain}, duality theory \cite{Ber_nlp_book:1999} shows that:
\begin{align}\label{eq_duality_theory}
 {q'_t}^* = {f'_t}^*
\end{align}
with ${q'_t}^* = \max_{\mu \in \rr^m_+} q'(\mu)$ and ${f'_t}^*$ the minimum of \eqref{eq_opt1}--\eqref{eq_opt1_domain}.

Thanks to this result, instead of minimizing the primal problem, we may maximize the dual problem, which is often more amenable to decomposition due to simpler constraints. Since we may not have the gradient of $q'$ in all points of $\rr^m_+$, we will use a method based on the subgradient.

\begin{definition}
A vector $d$ is called a \textit{subgradient} of a convex function $f$ over $\Xc$ at the point $x \in \Xc$ if:
\begin{align}
f(y) \geq f(x) + (y-x)^T d, \quad \forall y \in \Xc
\end{align}

The set of all subgradients of $f$ at the point $x$ is called the subdifferential of $f$ at $x$, denoted $\partial f(x)$.

\end{definition}

For each Lagrange multiplier $\bar{\mu} \in \rr^m_+$, first assume we have $\ub(\bar{\mu}) = \arg\min_{\ub \in \Omegab} \mathcal{L}'(\ub, \bar{\mu})$. Then a subgradient of the dual function is directly available, since \cite{Ber_nlp_book:1999}:
\begin{align}\label{eq_sub_dual_func}
q'(\mu) \leq q'(\bar{\mu}) + (\mu-\bar{\mu})^T g'(\ub(\bar{\mu})), \forall \mu \in \rr^m_+
\end{align}

In case an optimum of the Lagrangian is not attained due to termination of the optimization algorithm after a finite number of steps, a value $\tilde{\ub}(\bar{\mu})$ that satisfies
\begin{align}
 \mathcal{L}'(\tilde{\ub}(\bar{\mu}), \bar{\mu}) \leq \min_{\ub \in \Omegab} \mathcal{L}'(\ub, \bar{\mu}) + \delta
\end{align}
will lead to the following inequality:
\begin{align}\label{eq_approx_sub_dual_func}
q'(\mu) \leq q'(\bar{\mu}) + \delta + (\mu-\bar{\mu})^T g'(\tilde{\ub}(\bar{\mu})), \forall \mu \in \rr^m_+
\end{align}
where $g'(\tilde{\ub}(\bar{\mu}))$ is called $\delta$-subgradient of the dual function $q$ at the point $\bar{\mu}$. The set of all $\delta$-subgradients of $q$ at $\bar{\mu}$ is called $\delta$-subdifferential of $q$ at $\bar{\mu}$.

This means we do not have to look for a subgradient (or $\delta$-subgradient) of the dual function, it is available by just evaluating the constraint function at the primal value $\ub(\bar{\mu})$ (or $\tilde{\ub}(\bar{\mu})$).

\subsection{The main algorithm}
We organize our algorithm for solving \eqref{eq_opt0}--\eqref{eq_opt0_domain} at time step $t$ in a nested iteration of an outer and inner loop. The main procedure is described as follows:

\hrule
\medskip
\begin{algorithm}\label{alg_avg_grad_jacobi}
\textbf{Approximate subgradient method with nested Jacobi iterations}
\begin{enumerate}
 \item Given a Slater vector $\bar{\ub}_{t}$ of \eqref{eq_opt0}--\eqref{eq_opt0_domain}, determine $c_t$ and construct the tightened problem \eqref{eq_opt1}--\eqref{eq_opt1_domain}.
 \item Determine step size $\alpha_t$ and suboptimality $\varepsilon_t$, see later in Section~\ref{sec_outer_loop_params}.
 \item Determine $\bar{k}_t$ (the sufficient number of outer iterations), see later in Section~\ref{sec_outer_loop_kbar}.
 \item \textbf{Outer loop:} Set $\mu^{(0)} = 0\cdot\mathbf{1}_m$. For $k=0,\dots,\bar{k}$, find $\ub^{(k)}, \mu^{(k+1)}$ such that:
\begin{align}
\mathcal{L}'(\ub^{(k)},\mu^{(k)}) &\leq \min_{\ub \in \Omegab} \mathcal{L}'(\ub,\mu^{(k)}) + \varepsilon_t \label{eq_subgrad_iter1}\\ 
\mu^{(k+1)} &= \mathcal{P}_{\rr^m_+} \bigg\lbrace \mu^{(k)} + \alpha_t d^{(k)} \bigg\rbrace \label{eq_subgrad_iter2}
\end{align}
where $\mathcal{P}_{\rr^m_+}$ denotes the projection onto the nonnegative orthant, $d^{(k)} = g'\big(\ub^{(k)},x_t\big)$.

\textbf{Inner loop:}
 \begin{itemize}
  \item Determine $\bar{p}_k$ (the sufficient number of inner iterations), see later in Section~\ref{sec_inner_loop_pbar}.
  \item  Solve problem \eqref{eq_subgrad_iter1} in a distributed way with a Jacobi algorithm. For $p=0,\dots,\bar{p}_k$, every subsystem $i$ computes:
\begin{align}\label{eq_Jacobi_iter}
\ub^i(p+1) =& \arg\min_{\ub^i \in \Omegab_i} \mathcal{L}'(\ub^1(p),\dots,\ub^{i-1}(p),\ub^i,\nonumber\\
&\qquad \ub^{i+1}(p),\dots,\ub^{M}(p), \mu^{(k)})
\end{align}
where $\Omegab_i$ is the local constraint set for control variables of subsystem $i$.
  \item Define $\ub^{(k)} \triangleq [\ub^1(\bar{p}_k)^T,\dots,\ub^{M}(\bar{p}_k)^T]^T$, which is guaranteed to satisfy \eqref{eq_subgrad_iter1}.
 \end{itemize}

 \item Compute $\hat{\ub}^{(\bar{k}_t)} = \frac{1}{\bar{k}_t}\sum_{l=0}^{\bar{k}_t} \ub^{(l)}$, take $\ub_t = \hat{\ub}^{(\bar{k}_t)}$ as the solution of \eqref{eq_opt0}--\eqref{eq_opt0_domain}.

\end{enumerate}
\end{algorithm}
\hrule
\medskip

\begin{remark}
Algorithm \ref{alg_avg_grad_jacobi} is suitable for implementation in a hierarchical fashion where the main computations occur in the Jacobi iterations and are executed in parallel by local controllers, while the updates of dual variables and common parameters are carried out by a higher-level coordinating controller. In the inner loop, each subsystem only needs to communicate with its neighbors, which will be discussed in Section~\ref{sec_dist_Jacobi}. This algorithm is also amenable to implementation in distributed settings, where there are communication links available to help determine and propagate the common parameters $\alpha_t, \varepsilon_t, \bar{k}_t$, and $\bar{p}_k$.
\end{remark}

In the following sections, we will describe in detail how the computations are derived, and what the resulting properties are.

\subsection{Outer loop: Approximate subgradient method}

The outer loop at iteration $k$ uses an approximate subgradient method. The primal average sequence $\hat{\ub}^{(k)} = \frac{1}{k}\sum_{l=0}^k \ub^{(l)}$ has the following properties:

\begin{align}
\mathrm{For~} k \geq 1: \nonumber\\
 \left\lVert\Big[g'\Big(\hat{\ub}^{(k)},x_t\Big)\Big]^+\right\rVert_2 &\leq \frac{1}{k \alpha_t} \bigg(\frac{3}{\gamma_t}[f(\bar{\ub}_t,x_t)-{q'_t}^*] \nonumber\\
 &\qquad+ \frac{\alpha_t {L'_t}^2}{2\gamma_t} + \alpha_t {L'_t} \bigg) \label{eq_const_bound}\\
 f\Big(\hat{\ub}^{(k)},x_t\Big) &\leq {f'_t}^* + \frac{\big\|\mu^{(0)}\big\|^2_2}{2k\alpha_t} + \frac{\alpha_t {L'_t}^2}{2} + \varepsilon_t \label{eq_cost_bound}
\end{align}
where ${g'}^+$ denotes the constraint violation, i.e. ${g'}^+ = \max\{g',0\cdot \mathbf{1}_m\}$. The proof of \eqref{eq_const_bound} can be found in \cite{NedOzd_approx:2009}, and the proof of \eqref{eq_cost_bound} is given in Appendix~\ref{app_cost_bound}.

\subsubsection{Determining $\alpha_t$ and $\varepsilon_t$}\label{sec_outer_loop_params}

Using the lower bound of the cost reduction \eqref{eq_cost_decrease} and the upper bound of the suboptimality \eqref{eq_cost_bound} for the tightened problem \eqref{eq_opt1}--\eqref{eq_opt1_domain}, we will choose $\alpha_t$ and $\varepsilon_t$ such that $f(\ub_{t},x_t) < f(\ub_{t-1},x_{t-1})$.

The step size $\alpha_t$ and suboptimality $\varepsilon_t$ should satisfy:
\begin{align} \label{eq_alpha_epsilon_t}
 \frac{\alpha_t {L'_t}^2}{2} + \varepsilon_t \leq \Delta_t
\end{align}
where $\Delta_t$ is defined in \eqref{eq_def_Delta}, and ${L'_t}$ is the norm bound for $g'$.
This condition allows us to show the decreasing property of the cost function in problem \eqref{eq_cmpc}--\eqref{eq_cmpc_xt}, which can then be used as a Lyapunov function.

Note that a larger $\alpha_t$ will lead to a smaller number of outer iterations, while a larger $\varepsilon_t$ will lead to a smaller number of inner iterations. 
For the remainder of the paper we choose their values according to
\begin{align}
\alpha_t &= \frac{\Delta_t}{{L'_t}^2} \label{eq_alpha_t} \\
\varepsilon_t &= \frac{\Delta_t}{2} \label{eq_approx_subgrad_t}
\end{align}

\subsubsection{Determining $\bar{k}_t$}\label{sec_outer_loop_kbar}

Using the constraint violation bound \eqref{eq_const_bound}, we will choose $\bar{k}_t$ such that at the end of the algorithm, we will get a feasible solution for problem \eqref{eq_opt0}--\eqref{eq_opt0_domain}, which is the average of primal iterates generated by \eqref{eq_subgrad_iter1}:
\begin{align}\label{eq_primal_average}
\hat{\ub}^{(\bar{k}_t)} = \frac{1}{\bar{k}_t}\sum_{l=0}^{\bar{k}_t} \ub^{(l)}
\end{align}

The subgradient iteration \eqref{eq_subgrad_iter1}--\eqref{eq_subgrad_iter2} is performed for $k=1,\dots,\bar{k}_t$, with the integer
\begin{align}\label{eq_kbar}
 \bar{k}_t = \bigg\lceil \frac{1}{\alpha_t c_t} \bigg( &\frac{3}{\gamma_t} f(\bar{\ub}_t,x_t) + \frac{\alpha_t {L'_t}^2}{2 \gamma_t} + \alpha_t {L'_t} \bigg) \bigg\rceil
\end{align}
defined \emph{a priori}, where $\lceil\cdot\rceil$ is the ceiling operator which gives the closest integer equal to or above a real value, $\gamma_t = \min_{j=1,\dots,m}\{-g'_{j}(\bar{\ub}_t,x_t)\}= \min_{j=1,\dots,m}\{-g_{j}(\bar{\ub}_t,x_t)\} - c_t$, and $\bar{\ub}_t$ is the Slater vector of \eqref{eq_opt1}--\eqref{eq_opt1_domain}.


\subsection{Inner loop: Jacobi method}

The inner iteration \eqref{eq_Jacobi_iter} performs parallel local optimizations based on a standard Jacobi distributed optimization method for a convex function $\mathcal{L}'(\ub,\mu^{(k)})$ over a Cartesian product, as described in \cite[Section 3.3]{BerTsi_pdc:1989}. In order to find the sufficient stopping condition of this Jacobi iteration, we need to characterize the convergence rate of this algorithm. In the following, we summarize the condition for convergence of the Jacobi iteration, noting that $\mathcal{L}'(\ub,\mu^{(k)})$ is a strongly convex quadratic function with respect to $\ub$.

\begin{proposition}\label{prop_jacobi_convergence}
Suppose the following condition holds:
\begin{align} \label{eq_cond_contract_coupling}
\lambda_{\min}(H_{ii}) > \sum_{j \neq i} \bar{\sigma}(H_{ij}), \forall i
\end{align}
where $H_{ij}$ with $i,j \in \{1, \dots, M\}$ denotes a submatrix of the Hessian $H$ of $\mathcal{L}'$ w.r.t.~$\ub$, containing entries of $H$ in rows belonging to subsystem $i$ and columns belonging to subsystem $j$, $\lambda_{\min}$ means the smalleast eigenvalue, and $\bar{\sigma}$ denotes the maximum singular value.

Then $\exists \phi \in (0,1)$ such that the aggregate solution of the Jacobi iteration \eqref{eq_Jacobi_iter} satisfies:
\begin{align} \label{eq_Jacobi_converge_rate_norm_2}
 \|\ub(p) - \ub^*\|_{2} \leq M \phi^p \max_i \|\ub^i(0) - \ub^{i*}\|_{2}, \quad \forall p \geq 1
\end{align}
where $\ub^*=\arg\min_{\ub \in \Omegab} \Lc'(\ub, \mu^{(k)})$, and $\ub^{i*}$ is the component of subsystem $i$ in $\ub^*$.
\end{proposition}

We provide a proof for Proposition~\ref{prop_jacobi_convergence} in Appendix~\ref{app_proof_jacobi_convergence}.

\begin{remark}
This proposition provides a linear convergence rate of the Jacobi iteration, under the condition of \emph{weak dynamical couplings} between subsystems. For the sake of illustrating condition \eqref{eq_cond_contract_coupling}, let all subsystems have the same number of inputs. Consequently, $H_{ij}$ is a square and symmetric matrix for each pair $(i,j)$, hence the maximum singular value $\bar{\sigma}(H_{ij})$ equals to the maximum eigenvalue. Inequality \eqref{eq_cond_contract_coupling} thus reads:
\begin{align*}
\lambda_{\min}(H_{ii}) > \sum_{j \neq i} \lambda_{\max}(H_{ij}), \forall i
\end{align*}
which implies that the couplings represented by $H$ are small in comparison with each local cost.
\end{remark}

\begin{remark}
Note that the strong convexity of $\Lc'$ and the condition \eqref{eq_cond_contract_coupling} are required only for the convergence rate result of the Jacobi iteration in which $\Lc'$ is a quadratic function.
Extensions to other types of systems, where the Lagrangian can be solved with bounded suboptimality, are immediate.
In such cases we simply need to replace the Jacobi iteration with the new algorithm in the inner loop, while the outer loop will remain intact.
\end{remark}

\subsubsection{Determining $\bar{p}_k$}\label{sec_inner_loop_pbar}

As $\mathcal{L}'(\ub,\cdot)$ is continuously differentiable in a closed bounded set $\Omegab$, it is Lipschitz continuous.

Suppose we know the Lipschitz constant $\Lambda$ of $\mathcal{L}'(\ub,\cdot)$ over $\Omegab$, i.e. for any $\ub^1, \ub^2 \in \Omegab$ the following inequality holds:
\begin{align} \label{eq_Jacobi_Lipschitz_const}
  \|\mathcal{L}'(\ub^1,\mu^{(k)}) - \mathcal{L}'(\ub^2,\mu^{(k)})\|_2 \leq \Lambda \|\ub^1 - \ub^2\|_2
\end{align}

Taking $\ub^1 = \ub(\bar{p}_k)$ and $\ub^2 = \ub^*$ in \eqref{eq_Jacobi_Lipschitz_const}, and combining it with \eqref{eq_Jacobi_converge_rate_norm_2}, we obtain:
\begin{align} \label{eq_Jacobi_converge_cost}
 \|\mathcal{L}'(\ub(\bar{p}_k),\mu^{(k)}) - &\min_{\ub \in \Omegab} \mathcal{L}'(\ub,\mu^{(k)})\|_2 \leq \Lambda \|\ub(\bar{p}_k) - \ub^*\|_2 \nonumber\\
&\leq \Lambda M \phi^{\bar{p}_k} \max_i \|\ub^i(0) - \ub^{i*}\|_{2}
\end{align}

For each $i \in \{1,\dots,M\}$, let $D_i$ denote the diameter of the set $\Omegab_i$ w.r.t.~the Euclidean norm, so we have $\|\ub^i(0) - \ub^{i*}\|_2 \leq D_i$. Hence the relation \eqref{eq_Jacobi_converge_cost} can be further simplified as
\begin{align} \label{eq_Jacobi_converge_simplified}
 \mathcal{L}'(\ub(\bar{p}_k),\mu^{(k)}) \leq \min_{\ub \in \Omegab} \mathcal{L}'(\ub,\mu^{(k)}) + \Lambda M \phi^{\bar{p}_k} \max_i D_i
\end{align}

Based on \eqref{eq_Jacobi_converge_simplified}, in order to use $\ub(\bar{p}_k)$ as the solution $\ub^{(k)}$ that satisfies \eqref{eq_subgrad_iter1}, we choose the smallest integer $\bar{p}_k$ such that $\Lambda M \phi^{\bar{p}_k} \max_i D_i \leq \varepsilon_t$:
\begin{align} \label{eq_p_bar}
 \bar{p}_k = \bigg \lceil \log_\phi \frac{\varepsilon_t}{\Lambda M \max_i D_i} \bigg \rceil
\end{align}

\section{Properties of the algorithm} \label{sec_properties}

\subsection{Distributed Jacobi algorithm with guaranteed convergence}\label{sec_dist_Jacobi}

The computations in the inner loop can be executed by subsystems in parallel. Let us define an $r$-step extended neighborhood of a subsystem $i$, denoted by $\Nr^i_r$, as the set containing all subsystems that can influence subsystem $i$ within $r$ successive time steps. $\Nr^i_r$ is the union of subsystem indices in the neighborhoods of all subsystems in $\Nr^i_{r-1}$:
\begin{align}
\Nr^i_r = \bigcup_{j \in \Nr^i_{r-1}} \Nr^j
\end{align}
where $\Nr^i_1 = \Nr^i$. We can see that in order to get update information in the Jacobi iterations, each subsystem $i$ needs to communicate only with subsystems in $\Nr_{N-1}^i$, where $N$ is the prediction horizon. This set includes all other subsystems that couple with $i$ in the problem \eqref{eq_opt0}--\eqref{eq_opt0_domain} after eliminating the state variables. 
This communication requirement indicates that we will benefit from communication reduction when the number of subsystems $M$ is much larger than the horizon $N$, and the coupling structure is sparse.

Assume that the weak coupling condition \eqref{eq_cond_contract_coupling} holds, then after $\bar{p}_k$ iterations as computed by \eqref{eq_p_bar}, the Jacobi algorithm generates a solution $\ub^{(k)} \triangleq \ub(\bar{p}_k)$ that satisfies \eqref{eq_subgrad_iter1} in the outer loop.

\subsection{Feasible primal solution}
\begin{proposition}\label{prop_feas_primal}
 Suppose Assumptions \ref{asm_decentralized_controller} and \ref{asm_cond_slater} hold. Construct $g'$ as in \eqref{eq_tight_const}, $\alpha_t$ as in \eqref{eq_alpha_t}. Let the outer loop \eqref{eq_subgrad_iter1}--\eqref{eq_subgrad_iter2} with $\mu^{(0)}=0 \cdot \mathbf{1}_m$ be iterated for $k=0,\dots,\bar{k}_t$. Then $\hat{\ub}^{(\bar{k}_t)}$ is a feasible solution of \eqref{eq_opt0}--\eqref{eq_opt0_domain}, where $\hat{\ub}^{(\bar{k}_t)}$ is the primal average, computed by \eqref{eq_primal_average}.
\end{proposition}

\textbf{Proof}: With a finite number of $\bar{k}_t$ iterations \eqref{eq_const_bound} reads as
\begin{align}
\left\lVert\Big[g'\Big(\hat{\ub}^{(\bar{k}_t)},x_t\Big)\Big]^+\right\rVert_2 \leq \frac{1}{\bar{k}_t \alpha_t} &\bigg(\frac{3}{\gamma_t}\big[f(\bar{\ub}_t,x_t) - {q'_t}^*\big] \nonumber\\
 &+ \frac{\alpha_t {L'_t}^2}{2\gamma_t} + \alpha_t {L'_t} \bigg)
\end{align}

Moreover, the dual function $q'_t$ is a concave function, therefore ${q'_t}^* \geq q'(0,x_t)$. Recall that $f(\ub,x_t) > 0, \forall \ub \neq 0, x_t \neq 0$, thus $q'(0,x_t) = \min_{\ub \in \Omegab} f(\ub,x_t) + 0\cdot \mathbf{1}_m^T g'(\ub,x_t) = \min_{\ub \in \Omegab} f(\ub,x_t) > 0$, thus
\begin{align}\label{eq_const_bound_kbar}
 \left\lVert\Big[g'\Big(\hat{\ub}^{(\bar{k}_t)},x_t\Big)\Big]^+\right\rVert_2 < \frac{1}{\bar{k}_t \alpha_t} &\bigg(\frac{3}{\gamma_t}f(\bar{\ub}_t,x_t) \nonumber\\
 &+ \frac{\alpha_t {L'_t}^2}{2\gamma_t} + \alpha_t {L'_t} \bigg)
\end{align}

Combining \eqref{eq_const_bound_kbar} with \eqref{eq_kbar}, and noticing that $\bar{k}_t$ and $c_t$ are all positive lead to
\begin{align}
 \left\lVert\Big[g'\Big(\hat{\ub}^{(\bar{k}_t)},x_t\Big)\Big]^+\right\rVert_2 &< c_t \\
\Rightarrow g'_j\Big(\hat{\ub}^{(\bar{k}_t)},x_t\Big) &< c_t, \quad j=1,\dots,m \\
\Rightarrow g_j\Big(\hat{\ub}^{(\bar{k}_t)},x_t\Big) &< 0, \quad j=1,\dots,m \label{eq_const_strictly_feasible}
\end{align}
where the last inequality implies that $\hat{\ub}^{(\bar{k}_t)}$ is a feasible solution of problem \eqref{eq_opt0}--\eqref{eq_opt0_domain},
due to $c_t < \min_{j=1,\dots,m}\{-g_j(\bar{\ub}_t,x_t)\}$.\cvd

\subsection{Closed-loop stability}
\begin{proposition}\label{prop_stability}
 Suppose Assumptions \ref{asm_cond_slater}, \ref{asm_cost_decrease}, and  \ref{asm_g_norm_bound} hold. Then the solution $\hat{\ub}^{(\bar{k}_t)}$ generated by Algorithm~\ref{alg_avg_grad_jacobi} satisfies the following inequality:
\begin{align}
f(\ub_t,x_t) < f(\ub_{t-1},x_{t-1}), \quad \forall t \in \zz_+
\end{align}
\end{proposition}

\textbf{Proof:} Using \eqref{eq_cost_bound} and \eqref{eq_alpha_epsilon_t}, and noting that $\mu^{(0)}=0$, we obtain:
\begin{align} \label{eq_cost_ineq1}
f\Big(\hat{\ub}^{(\bar{k}_t)},x_t\Big) \leq {f'_t}^* + \frac{\|\mu^{(0)}\|}{2\bar{k}_t\alpha_t} + \frac{\alpha_t {L'_t}^2}{2} + \varepsilon_t \leq {f'_t}^* + \Delta_t
\end{align}

Notice that $\bar{\ub}_t$ is also a feasible solution of \eqref{eq_opt1}--\eqref{eq_opt1_domain} (due to the way we construct the tightened problem: $\bar{\ub}_t$ still belongs to the interior of the tightened constraint set), while ${f'_t}^*$ is the optimal cost value of this problem. As a consequence,
\begin{align} \label{eq_cost_ineq2}
{f'_t}^*  \leq f(\bar{\ub}_t,x_t)
\end{align}

Combining \eqref{eq_cost_ineq1}, \eqref{eq_cost_ineq2}, and \eqref{eq_cost_decrease}, and noting that $\ub_{t} = \hat{\ub}^{(\bar{k}_t)}$ leads to:
\begin{align}
f(\ub_{t},x_t) < f(\ub_{t-1},x_{t-1}), \quad \forall t \in \zz_+
\end{align}\cvd

Note that besides the decreasing property of $f(\ub_{t},x_t)$, all the other conditions for Lyapunov stability of MPC \cite{MayRaw:00} are satisfied. Therefore, Proposition~\ref{prop_stability} leads to closed-loop MPC stability, where the cost function $f(\ub_{t},x_t)$ is a Lyapunov candidate function.

\section{Realization of the assumptions}

In this section, we discuss the method to update the Slater vector and the constraint norm bound for each time step, implying that Assumptions~\ref{asm_cond_slater} and \ref{asm_g_norm_bound} are only necessary in the first time step ($t=0$).

\subsection{Updating the Slater vector}\label{sec_slate_update}
\begin{lemma}
 Suppose Assumption \ref{asm_pos_inv_terminal_set} holds. Let $\ub_t$ be the solution of the MPC problem \eqref{eq_cmpc}--\eqref{eq_cmpc_xt} at time step $t$, computed by Algorithm~\ref{alg_avg_grad_jacobi}. Then $\tilde{\ub}_{t+1}$ constructed by shifting $\ub_{t}$ one step ahead and adding $\tilde{u}_{t+N}=Kx_{t+N}$, is a Slater vector for constraint \eqref{eq_opt0_const} at time step $t+1$.
\end{lemma}

\textbf{Proof:} Note that based on Proposition~\ref{prop_feas_primal}, $\hat{\ub}^{(\bar{k}_t)}$ is a feasible solution of problem \eqref{eq_opt0}--\eqref{eq_opt0_domain}. Moreover, the strict inequality \eqref{eq_const_strictly_feasible} means that $\hat{\ub}^{(\bar{k}_t)}$ is in the interior of the constraint set of \eqref{eq_cmpc}--\eqref{eq_cmpc_xt}. This also yields:
\begin{align}
 x_{t+N} \in \mathrm{int}(\Xf)
\end{align}

Moreover, due to Assumption \ref{asm_pos_inv_terminal_set}, we have $(A + BK)x_{t+N} \in \mathrm{int}(\Xf)$. This means that if we use $\tilde{u}_{t+N}=Kx_{t+N}$, then the next state is also in the interior of the terminal constraint set $\Xf$. Note that $\Uc$ and $\Xc$ do not change when problem \eqref{eq_cmpc}--\eqref{eq_cmpc_xt} is shifted from $t$ to $t+1$, hence all the inputs of $\tilde{\ub}_{t+1}$ and their subsequent states are in the interior of the corresponding constraint sets. Therefore, $\tilde{\ub}_{t+1}$ as constructed at step 5 of Algorithm~\ref{alg_avg_grad_jacobi} is a Slater vector for the constraint \eqref{eq_opt0_const} at time step $t+1$.\cvd

This means we can use $\bar{\ub}_{t+1} = \tilde{\ub}_{t+1}$ as the qualifying Slater vector for Assumption~\ref{asm_cond_slater} at time step $t+1$.

\subsection{Updating the constraint norm bound}\label{sec_norm_bound_update}

In our general problem setup, $g(\ub, x)$ is composed of affine functions over $\ub$ and $x$, and thus can be written compactly as
\begin{align}
 g(\ub, x) = \Xi x + \Theta \ub + \tau
\end{align}
with constant matrices $\Xi, \Theta$ and vector $\tau$. Then for each $x_{t-1}$, $x_t$, and $\ub \in \Omegab$,
the following holds:
\begin{align}
 g(\ub, x_t) = g(\ub, x_{t-1}) + \Xi (x_t - x_{t-1}) \nonumber \\
 \Rightarrow\|g(\ub, x_t)\|_2 \leq \|g(\ub, x_{t-1})\|_2 + \|\Xi (x_t - x_{t-1})\|_2
\end{align}

In order to find a bound $L_t$ for $g(\ub, x_t)$ in each $t\geq1$ step, we assume to have the constraint norm bound available
from the previous step:
\begin{align}
 L_{t-1} \geq \|g(\ub, x_{t-1})\|_2, \forall \ub \in \Omegab
\end{align}

Hence, combining the above inequalities a norm bound update for $g(\ub, x_t)$ can be obtained as:
\begin{align}
L_t = L_{t-1} + \|\Xi (x_t - x_{t-1})\|_2
\end{align}

\section{Conclusions}\label{sec_conclusion}

We have presented a constraint tightening approach for solving an MPC optimization problem with guaranteed feasibility and stability after a finite number of iterations. The new method is applicable to large-scale systems with coupling in dynamics and constraints, and the solution is based on approximate subgradient and Jacobi iterative methods, which facilitate implementation in a hierarchical or distributed way. Future extensions of this scheme include \emph{a posteriori} choice of the solution by comparing the cost functions associated with the Slater vector $\bar{\ub}_t$ and the primal average $\hat{\ub}^{(\bar{k}_t)}$ in a distributed way.


\section*{Acknowledgement}

The authors would like to thank Ion Necoara for helpful discussions on the topic of this paper.

Research supported by the European Union Seventh Framework STREP project ``Hierarchical and distributed model predictive control (HD-MPC)", contract number INFSO-ICT-223854, and the European Union Seventh Framework Programme [FP7/2007-2013] under grant agreement no. 257462 HYCON2 Network of Excellence.


\bibliographystyle{ieeetr}
\bibliography{Doan_bib}

\begin{thebibliography}{10}

\bibitem{rawlings_coordinating_2008}
J.~B. Rawlings and B.~T. Stewart, ``Coordinating multiple optimization-based
  controllers: New opportunities and challenges,'' {\em Journal of Process
  Control}, vol.~18, pp.~839--845, Oct. 2008.

\bibitem{Mac:02}
J.~M. Maciejowski, {\em Predictive Control with Constraints}.
\newblock Harlow, England: Prentice-Hall, 2002.

\bibitem{CamBor:99}
E.~F. Camacho and C.~Bordons, {\em Model Predictive Control}.
\newblock London: Springer, 1999.

\bibitem{RawMay:09}
J.~B. Rawlings and D.~Q. Mayne, {\em Model Predictive Control: Theory and
  Design}.
\newblock Madison, WI: Nob Hill Publishing, 2009.

\bibitem{scattolini_architecture_2009}
R.~Scattolini, ``Architectures for distributed and hierarchical model
  predictive control - {A} review,'' {\em Journal of Process Control}, vol.~19,
  pp.~723--731, May 2009.

\bibitem{Venkat_tcst:2008}
A.~Venkat, I.~Hiskens, J.~Rawlings, and S.~Wright, ``Distributed {MPC}
  strategies with application to power system automatic generation control,''
  {\em IEEE Transactions on Control Systems Technology}, vol.~16,
  pp.~1192--1206, Nov. 2008.

\bibitem{yuji_wakasa_decentralized_2008}
Y.~Wakasa, M.~Arakawa, K.~Tanaka, and T.~Akashi, ``Decentralized model
  predictive control via dual decomposition,'' in {\em 47th IEEE Conference on
  Decision and Control}, pp.~381--386, 2008.

\bibitem{NecSuy_prox:2008}
I.~Necoara and J.~Suykens, ``Application of a smoothing technique to
  decomposition in convex optimization,'' {\em IEEE Transactions on Automatic
  Control}, vol.~53, pp.~2674--2679, Dec. 2008.

\bibitem{DoaKev:09}
D.~Doan, T.~Keviczky, I.~Necoara, M.~Diehl, and B.~{De Schutter}, ``A
  distributed version of {Han}'s method for {DMPC} using local communications
  only,'' {\em Control Engineering and Applied Informatics}, vol.~11,
  pp.~6--15, Sept. 2009.

\bibitem{DoaKev:11IFACWC}
M.~D. Doan, T.~Keviczky, and B.~{De Schutter}, ``A dual decomposition-based
  optimization method with guaranteed primal feasibility for hierarchical {MPC}
  problems,'' in {\em 18th IFAC World Congress}, (Milan, Italy), Aug. 2011.

\bibitem{Kuwata:07}
Y.~Kuwata, A.~Richards, T.~Schouwenaars, and J.~P. How, ``Distributed robust
  receding horizon control for multivehicle guidance,'' {\em IEEE Transactions
  on Control Systems Technology}, vol.~15, 2007.

\bibitem{Ber_nlp_book:1999}
D.~P. Bertsekas, {\em Nonlinear programming}.
\newblock Belmont, MA: Athena Scientific, 1999.

\bibitem{Siljak:78_lss}
D.~D. \v{S}iljak, {\em Large-scale dynamic systems: Stability and structure}.
\newblock New York, NY: North Holland, 1978.

\bibitem{NedOzd_approx:2009}
A.~Nedic and A.~Ozdaglar, ``Approximate primal solutions and rate analysis for
  dual subgradient methods,'' {\em SIAM Journal on Optimization}, vol.~19,
  pp.~1757--1780, Nov. 2009.

\bibitem{BerTsi_pdc:1989}
D.~P. Bertsekas and J.~N. Tsitsiklis, {\em Parallel and Distributed
  Computation: Numerical Methods}.
\newblock Upper Saddle River, NJ: Prentice-Hall, 1989.

\bibitem{MayRaw:00}
D.~Q. Mayne, J.~B. Rawlings, C.~V. Rao, and P.~O.~M. Scokaert, ``Constrained
  model predictive control: Stability and optimality,'' {\em Automatica},
  vol.~36, pp.~789--814, June 2000.

\end{thebibliography}

\section{Appendix}

\subsection{Proof of the upper bound on the cost function \eqref{eq_cost_bound}} \label{app_cost_bound}

This proof is an extension of the proof of Proposition~3(b) in \cite{NedOzd_approx:2009}, the main difference being the incorporation of the suboptimality $\varepsilon_t$ in the update of the primal variable \eqref{eq_subgrad_iter1}.

Using the convexity of the cost function, we have:
\begin{align}
 f(\hat{\ub}^{(k)}) = f\Bigg(\frac{1}{k}\sum_{l=0}^{k-1} \ub^{(l)}\Bigg) \leq \frac{1}{k}\sum_{l=0}^{k-1} f(\ub^{(l)}) \nonumber\\
= \frac{1}{k}\sum_{l=0}^{k-1} \big(f(\ub^{(l)}) + (\mu^{(l)})^T g'(\ub^{(l)}) \big) - \frac{1}{k}\sum_{l=0}^{k-1}(\mu^{(l)})^T g'(\ub^{(l)})
\end{align}

Note that $\mathcal{L}'\big(\ub^{(l)},\mu^{(l)}\big) = \bigg(f(\ub^{(l)}) + g'(\ub^{(l)})^T\mu^{(l)}\bigg)$ and
\begin{align}
\mathcal{L}'\big(\ub^{(l)},\mu^{(l)}\big) \leq \min_{\ub \in \Omegab} \mathcal{L}'\big(\ub^{(l)},\mu^{(l)}\big) + \varepsilon_t = q'\big(\mu^{(l)}\big) + \varepsilon_t,\nonumber\\
 \forall l < k
\end{align}

Combining the two inequalities above, we then have:
\begin{align} \label{eq_cost_bound_proof1}
f(\hat{\ub}^{(k)}) &\leq \frac{1}{k}\sum_{l=0}^{k-1} q'\big(\mu^{(l)}\big) + \varepsilon_t - \frac{1}{k}\sum_{l=0}^{k-1}(\mu^{(l)})^T g'(\ub^{(l)}) \nonumber \\
&\leq {q'_t}^* + \varepsilon_t - \frac{1}{k}\sum_{l=0}^{k-1}(\mu^{(l)})^T d^{(l)}
\end{align}
where $d^{(l)}=g'(\ub^{(l)})$, and the last inequality is due to ${q'_t}^* \geq q'\big(\mu^{(l)}\big), \forall l$.

Using the expression of squared sum:
\begin{align}
 \|\mu^{(l+1)}\|_2^2 &\leq \|\mu^{(l)} + \alpha_t d^{(l)}\|_2^2 \nonumber\\
 &= \|\mu^{(l)}\|_2^2 + 2\alpha_t (\mu^{(l)})^T d^{(l)} + \|\alpha_t d^{(l)}\|_2^2
\end{align}
we have:
\begin{align}
- (\mu^{(l)})^T d^{(l)}  \leq \frac{1}{2\alpha_t} \bigg(\|\mu^{(l)}\|_2^2 - \|\mu^{(l+1)}\|_2^2 + \alpha_t^2\|d^{(l)}\|_2^2 \bigg)
\end{align}
for $l=0,\dots,k-1$.

Summing side by side for $l=0,\dots,k-1$, we get:
\begin{align} \label{eq_cost_bound_proof2}
 - \sum_{l=0}^{k-1}(\mu^{(l)})^T d^{(l)} \leq &\frac{1}{2\alpha_t} \bigg(\|\mu^{(0)}\|_2^2 - \|\mu^{(k)}\|_2^2\bigg) \nonumber\\
 &+ \frac{\alpha_t}{2}\sum_{l=0}^{k-1}\|d^{(l)}\|_2^2
\end{align}

Linking \eqref{eq_cost_bound_proof1} and \eqref{eq_cost_bound_proof2}, we then have:
\begin{align}
f(\hat{\ub}^{(k)}) &\leq {q'_t}^* + \varepsilon_t + \frac{1}{2k\alpha_t} \bigg(\|\mu^{(0)}\|^2 - \|\mu^{(k)}\|^2\bigg) \nonumber\\
& \quad + \frac{\alpha_t}{2k}\sum_{l=0}^{k-1}\|d^{(l)}\|^2 \nonumber \\
&\leq {q'_t}^* + \frac{\|\mu^{(0)}\|^2 }{2k\alpha_t} + \frac{\alpha_t {L'_t}^2}{2} + \varepsilon_t
\end{align}
in which we get the last inequality by using $L'_t$ as the norm bound for all $g'(\ub^{(l)}), l=0,\dots,k-1$.

Finally, with the Slater condition, there is no primal-dual gap, i.e. ${q'_t}^* = f_t^*$ (cf. \eqref{eq_duality_theory}), hence:
\begin{align*}
f(\hat{\ub}^{(k)}) \leq {f'_t}^* + \frac{\|\mu^{(0)}\|^2 }{2k\alpha_t} + \frac{\alpha_t {L'_t}^2}{2} + \varepsilon_t
\end{align*}
\cvd

\subsection{Proof of the convergence result of the Jacobi iteration (Proposition~\ref{prop_jacobi_convergence})}\label{app_proof_jacobi_convergence}
 According to Proposition 3.10 in \cite[Chapter 3]{BerTsi_pdc:1989}, the Jacobi algorithm has a linear convergence w.r.t. the block-maximum norm, as defined below:

\begin{definition}
 For each vector $x = [x_1^T, \dots, x_M^T]$ with $x_i \in \rr^{n_i}$, given a norm $\| \cdot \|_i$ for each $i$, the \textit{block-maximum norm} based on $\| \cdot \|_i$ is defined as:
\begin{align}
 \|x\|_{\textrm{b-m}} = \max_i \|x_i\|_i
\end{align}
\end{definition}

\begin{definition}
With any matrix $A \in \rr^{n_i \times n_j}$, we associate the \textit{induced matrix norm of the block-maximum norm}:
\begin{align}\label{eq_mat_blk_max_norm}
 \lVert A \rVert_{ij} = \max_{x \neq 0} \frac{\lVert A x \rVert_i}{\lVert x \rVert_j} = \max_{\lVert x \rVert_j = 1} \lVert A x \rVert_i
\end{align}
\end{definition}

In this paper, we use the Euclidean norm as the default basis for block-maximum norm, i.e. $\| \cdot \|_i = \| \cdot \|_2, \forall i$.

Proposition 3.10 in \cite[Chapter 3]{BerTsi_pdc:1989} states that $\ub(p)$ generated by \eqref{eq_Jacobi_iter} will converge to the optimizer of $\Lc'(\ub,x_t)$ with linear convergence rate w.r.t. block-maximum norm (i.e. $\|\ub(p) - \ub^*\|_{\textrm{b-m}} \leq \phi^{p}\|\ub(0) - \ub^*\|_{\textrm{b-m}}$, with $\ub^* = \arg\min_{\ub} \Lc'(\ub,x_t)$ and $\phi \in [0,1)$) if there exists a positive scalar $\gamma$ such that the mapping $R: \Omegab \mapsto \rr^{n_{\ub}}$, defined by $R(\ub) = \ub - \gamma \grad_{\ub} \Lc'(\ub,x_t)$, is a contraction w.r.t. the block-maximum norm.

Our focus now is to derive the condition such that $R(\ub)$ is a contraction mapping.

Note that since $f(\ub,x_t)$ is a quadratic function, and $g'(\ub,x_t)$ contains only linear functions, the function $\Lc'(\ub,x_t)$ is also a quadratic function w.r.t. $\ub$, hence it can be written as:

\begin{align}
 \Lc'(\ub,x_t) = \ub^T H \ub + b^T \ub + c
\end{align}
where $H$ is a symmetric, positive definite matrix, $b$ is a constant vector and $c$ is a constant scalar.

In order to derive the condition for $R(\ub)$ to be a contraction mapping, we will make use of Proposition 1.10 in \cite[Chapter 3]{BerTsi_pdc:1989}, stating that:

If $f: \rr^{n_{\ub}} \mapsto \rr^{n_{\ub}}$ is continuously differentiable and there exists a scalar $\phi \in [0,1)$ such that
\begin{align} \label{eq_cond_contract_mapping}
 \lVert I - \gamma G_i^{-1}\big(\grad_i F_i(\ub)\big)^T \rVert_{ii} + \sum_{j \neq i} \lVert \gamma G_i^{-1} \big( \grad_j F_i(\ub) \big)^T \rVert_{ij} \leq \phi, \nonumber\\
  \forall \ub \in \Omegab, \forall i
\end{align}
then the mapping $T: \Omegab \mapsto \rr^{n_{\ub}}$ defined with each component $i \in \{1, \dots, M\}$ by $T_i(\ub) = \ub_i - \gamma G_i^{-1} F(\ub)$ is a contraction with respect to the block-maximum norm.

The mapping $T(\ub)$ will become the mapping $R(\ub)$ if we take $G_i = I^{n_{\ub^i}}, \forall i$ and $F(\ub) = \grad_{\ub} \Lc'(\ub,x_t) = 2H \ub + b$. With such choice, and evaluating the induced matrix norm \eqref{eq_mat_blk_max_norm} in \eqref{eq_cond_contract_mapping}, the condition for contraction mapping of $R(\ub)$ is to find $\phi \in [0,1)$ such that:

\begin{align} \label{eq_cond_contract_norm_2}
 \lVert I^{n_{\ub^i}} - 2\gamma H_{ii} \rVert_2 + \sum_{j \neq i} \lVert 2\gamma H_{ij} \rVert_2 \leq \phi, \forall i
\end{align}
where $H_{ij}$ with $i,j \in \{1, \dots, M\}$ denotes the submatrix of $H$, containing entries at rows belonging to subsystem $i$ and columns belonging to subsystem $j$. Note that the matrix inside the first induced matrix norm is a square, symmetric matrix, while the matrices $H_{ij}$ are generally not symmetric, depending on the number of variables of each subsystem. The scalar $\phi \in [0,1)$ is also the modulus of the contraction.

Using the properties of eigenvalue and singular value of matrices, we transform \eqref{eq_cond_contract_norm_2} into the following inequality:

\begin{align} \label{eq_cond_contract_eig}
 \max_{\lambda} \lvert 2\gamma\lambda(H_{ii}) - 1 \rvert + 2\gamma \sum_{j \neq i} \bar{\sigma}(H_{ij}) \leq \phi, \forall i
\end{align}
where $\lambda$ means eigenvalue, and $\bar{\sigma}$ denotes the maximum singular value.

In order to find $\gamma > 0$ and $\phi \in [0,1)$ satisfying \eqref{eq_cond_contract_eig}, we need:

\begin{align}
 &\max_{\lambda} \lvert 2\gamma\lambda(H_{ii}) - 1 \rvert + 2\gamma \sum_{j \neq i} \bar{\sigma}(H_{ij}) < 1, \forall i \\
\Leftrightarrow &\left\lbrace \begin{array}{l}
                         2\gamma\lambda_{\max}(H_{ii}) - 1 + 2\gamma \sum_{j \neq i} \bar{\sigma}(H_{ij}) < 1 \\
                         1 - 2\gamma\lambda_{\min}(H_{ii}) + 2\gamma \sum_{j \neq i} \bar{\sigma}(H_{ij}) < 1
                        \end{array}\right., \forall i \\
\Leftrightarrow &\left\lbrace \begin{array}{l}
                         \gamma < 1/ \left(\lambda_{\max}(H_{ii}) + \sum_{j \neq i} \bar{\sigma}(H_{ij})\right) \\
                         \lambda_{\min}(H_{ii}) > \sum_{j \neq i} \bar{\sigma}(H_{ij})
                        \end{array}\right., \forall i \label{eq_cond_contract_final}
\end{align}

The first inequality of \eqref{eq_cond_contract_final} shows how to choose $\gamma$, while the second inequality of \eqref{eq_cond_contract_final} needs to be satisfied by the problem structure, which implies there are \emph{weak dynamical couplings} between subsystems.

In summary, the mapping $R(\ub)$ satisfies \eqref{eq_cond_contract_mapping} and thus is a contraction mapping if the following conditions hold:

\begin{enumerate}
 \item For all $i$:
\begin{align} \label{eq_cond_contract_coupling_appendix}
\lambda_{\min}(H_{ii}) > \sum_{j \neq i} \bar{\sigma}(H_{ij})
\end{align}
 \item The coefficient $\gamma$ is chosen such that:
\begin{align}\label{eq_cond_contract_gamma}
\gamma < \frac{1}{\lambda_{\max}(H_{ii}) + \sum_{j \neq i} \bar{\sigma}(H_{ij})}, \forall i
\end{align}
\end{enumerate}

So, when condition \eqref{eq_cond_contract_coupling_appendix} is satisfied and with $\gamma$ chosen by \eqref{eq_cond_contract_gamma}, we can define $\phi \in (0,1)$ as:
\begin{align}
\phi = \max_{i} \Bigg\lbrace \max \bigg\lbrace &2\gamma\big(\lambda_{\max}(H_{ii}) + \sum_{j \neq i} \bar{\sigma}(H_{ij})\big) - 1, \nonumber \\
&1 - 2\gamma \big(\lambda_{\min}(H_{ii}) - \sum_{j \neq i} \bar{\sigma}(H_{ij})\big) \bigg\rbrace \Bigg\rbrace
\end{align}

This $\phi$ is the modulus of the contraction $R(\ub)$, and also acts as the coefficient of the linear convergence rate of the Jacobi iteration \eqref{eq_Jacobi_iter}, which means:

\begin{align} \label{eq_Jacobi_converge_rate_bm_appendix}
 \|\ub(p) - \ub^*\|_{\textrm{b-m}} \leq \phi^p \|\ub(0) - \ub^*\|_{\textrm{b-m}}, \quad \forall p \geq 1
\end{align}
where $\ub^* = \arg\min_{\ub \in \Omegab} \mathcal{L}'(\ub, x_t)$.

Note that the closer of $\phi$ to 0, the faster the aggregate update $\ub(p)$ converges to the optimizer of the Lagrange function.

In order to get the convergence rate w.r.t. the Euclidean norm, we will need to link from the Euclidean norm to the block-maximum norm:
\begin{align}
 \|x\|_2 \leq \sum_{i=1}^M \|x^i\|_2 &\leq M \max_i \|x^i\|_2 = M \|x\|_{\textrm{b-m}}
\end{align}

Hence, the convergence rate of Jacobi iteration \eqref{eq_Jacobi_iter} w.r.t. the Euclidean norm is:

\begin{align} \label{eq_Jacobi_converge_rate_norm_2_appendix}
 \|\ub(p) - \ub^*\|_{2} \leq M \phi^p \max_i \|\ub^i(0) - \ub^{i*}\|_{2}, \quad \forall p \geq 1
\end{align}
\cvd

\end{document}